\documentclass[12pt]{article}
\usepackage{amsmath}
\usepackage{amssymb}
\usepackage{amscd}
\usepackage{amsthm}

\theoremstyle{plain}
\newtheorem{Thm}{Theorem}[section]
\newtheorem{Conj}[Thm]{Conjecture}
\newtheorem{Prop}[Thm]{Proposition}
\newtheorem{Cor}[Thm]{Corollary}
\newtheorem{Lem}[Thm]{Lemma}

\theoremstyle{definition}
\newtheorem{Defn}[Thm]{Definition}
\newtheorem{Expl}[Thm]{Example}

\newtheorem{Rem}[Thm]{Remark}

\numberwithin{equation}{section}

\title{Log Crepant Birational Maps and Derived Categories}
\author{Yujiro Kawamata}

\begin{document}

%\begin{abstract}
%We extend the conjecture on the derived equivalence for varieties with 
%equivalent canonical divisors to the logarithmic case and 
%prove it in the toric case. 
%\end{abstract}

%\subjclass[2000]{Primary 14E30, 18E30}

\maketitle

\section{Introduction}

The purpose of this paper is to extend the conjecture stated in the paper 
\cite{DK} to the logarithmic case and prove some supporting evidences.
\cite{DK}~Conjecture~1.2 predicts that birationally equivalent smooth 
projective varieties have equivalent derived categories if and only if 
they have equivalent canonical divisors.

According to the experience of the minimal model theory, one has to 
deal with singular varieties instead of only smooth varieties for the 
classification of algebraic varieties.  Moreover, we should consider not 
only varieties but also pairs consisting of varieties and divisors on them. 
These pairs are expected to have some milder singularities, the log terminal 
singularities.

On the other hand, the theory of derived categories works well under the 
smoothness assumption of the variety.
The reason is that the global dimension is finite only in the case of 
smooth varieties.

In this paper we shall consider pairs of varieties and $\mathbb{Q}$-divisors 
which have smooth local coverings (Definition~\ref{stack}), 
and conjecture that, 
if there is an equivalence of log canonical divisors between birationally 
equivalent pairs, then there is an equivalence of derived categories 
(Conjecture~\ref{conj}).
We note that we need to consider the sheaves on the stacks associated 
to the pairs instead of the usual sheaves on the varieties in order to 
have equivalences of derived categories 
as already notices in \cite{Francia}.
This is a generalization of the conjecture in \cite{DK}, and includes 
the case considered in \cite{BKR}.
We note that crepant resolutions in higher dimensions are rare but there are 
many log crepant partial resolutions of quotient singularities by finite
subgroups of general linear groups.
We note also that, even in the case in which there is a birational morphism 
between varieties, 
the direction of the inclusion of the category may be different from that
of the morphism, but coincides with that of the inequality of the log 
canonical divisors.

In \S 3, we consider the problem on recovering the variety from the category.
We prove that some basic birational invariants related to the 
canonical divisors can be recovered from the 
derived category (Theorem~\ref{derived}).
In particular, we prove the converse statement of the conjecture.
On the other hand, we remark that the variety itself can be reconstructed 
from the category of coherent sheaves (Theorem~\ref{coh}).

In \S 4, we consider the toroidal varieties and prove that the conjecture holds
in this case (Theorem~\ref{toroidal}).
This is a generalization of \cite{Francia}~Theorem~5.2.
We note that our result implies the McKay correspondence 
for abelian quotient singularities as a special case.

We conclude the paper with a remark on the relationship with the 
non-commutative geometry in Proposition~\ref{non-comm}.
We maybe need the moduli theoretic interpretation of the log crepant maps
in order to deal with the conjecture in the difficult general case.

The author would like to thank Yong-bin Ruan, Tom Bridgeland, 
Andrei Caldararu and Kenji Matsuki for useful discussions.

\section{Derived equivalence conjecture}

We shall consider pairs of varieties with $\mathbb{Q}$-divisors which 
have local coverings by smooth varieties:

\begin{Defn}\label{stack}
Let $X$ be a normal variety, and $B$ an effective $\mathbb{Q}$-divisor
on $X$ whose coefficients belong to the standard set 
$\{1 - 1/n ; n \in \mathbb{N}\}$.
Assume the following condition.

(*) There exists a quasi-finite and surjective morphism 
$\pi: U \to X$ from a smooth variety, which may be reducible, 
such that $\pi^*(K_X + B) = K_U$.
  
Let $R=(U \times_X U)\tilde{}$ be the normalization of the 
fiber product.
Then the projections $p_i: R \to U$ are etale for $i=1,2$, because there is
no ramification divisor for $p_i$.

We define the associated Deligne-Mumford stack $\mathcal{X}$ 
as a $2$-functor 
\[
\mathcal{X}: (Sch) \to (Groupoid)
\]
where $(Sch)$ is the category of schemes and 
$(Groupoid)$ is the $2$-category of groupoids, categories themselves 
whose morphisms are only isomorphisms.
The $2$-functor $\mathcal{X}$ is the sheaf with respect to the etale topology 
which is defined in the following way.
For any scheme $B$, an object of the category $\mathcal{X}(B)$ 
is an element of the set $U(B) = \text{Hom}(B,U)$, and a morphism 
of the category $\mathcal{X}(B)$ is an element of the set $R(B) = 
\text{Hom}(B,R)$.
For a morphism $f: B' \to B$, we have a functor
$f^*: \mathcal{X}(B) \to \mathcal{X}(B')$ given by
$f^*: U(B) \to U(B')$ and $f^*: R(B) \to R(B')$.
\end{Defn}

If we replace $U$ by its etale covering, then the category 
$\mathcal{X}(B)$ is replaced by an equivalent category so that the stack 
$\mathcal{X}$ does not depend on the choice of the covering 
$\pi: U \to X$ but only on the pair $(X, B)$.

\begin{Conj}\label{conj}
Let $(X, B)$ and $(Y, C)$ be pairs of quasi-projective varieties 
with $\mathbb{Q}$-divisors which satisfy the condition (*) in 
Definition~\ref{stack}, and let $\mathcal{X}$ and $\mathcal{Y}$ be 
the associated stacks. 
Assume that there are proper
birational morphisms $\mu: W \to X$ and $\nu: W \to Y$ 
from a third variety $W$ such that $\mu^*(K_X+B) = \nu^*(K_Y+C)$.
Then there exists an equivalence as triangulated categories
$D^b(\text{Coh}(\mathcal{X})) \to D^b(\text{Coh}(\mathcal{Y}))$.
\end{Conj}

The pairs considered in the conjecture are very special kind of 
log terminal pairs.
But our assumption is sifficiently general in dimension $2$:

\begin{Prop}
(1) Let $(X, B)$ be a pair which satisfies the condition (*) in 
Definition~\ref{stack}. 
Then the pair $(X, B)$ is log terminal.

(2) Let $X$ be a normal surface, 
and $B$ an effective $\mathbb{Q}$-divisor 
on $X$ whose coefficients belong to the standard set 
$\{1 - 1/n ; n \in \mathbb{N}\}$.
Assume that the pair $(X, B)$ is log terminal.
Then the pair satisfies the condition (*) in Definition~\ref{stack}.  
\end{Prop}

\begin{proof}
(1) Let $\mu: Y \to X$ be a proper birational morphism, and $E$ 
an exceptional prime divisor of $\mu$.
Then there is a birational morphism $\nu: V \to U$ with a quasi-finite 
morphism $\rho: V \to Y$ and a exceptional prime divisor $F$ of $\rho$
which dominates $E$.
We can write $\mu^*(K_X+B)=K_Y+aE+\dots$ and 
$\nu^*K_U=K_V+bF+\dots$.
If $e$ is the ramification index of the morphism $\rho$ along $F$,
then we have $b=ae+(e-1)$.
Since $b > 0$, we conclude that $a > -1$.

(2) We may always replace $X$ by its open covering in the etale topology
in the course of the proof.
Fixing a point $x \in X$, we shall 
construct our covering in a neighborhood of $x$.
Let $\{n_1, n_2, \dots \}$ be the set of all integers such that 
the numbers $1 - 1/n_i$ appear as coefficients of irreducible 
components of $B$ which pass through $x$.
We may assume that $n_1 > n_2 > \dots$.

By the classification of log terminal singularities on surfaces, 
$X$ has only quotient singularities in the etale topology.
Thus there is a quasi-finite and surjective morphism 
$\pi_1: X_1 \to X$ from a smooth surface which is etale in codimension $1$.
Let $B_1 = \pi_1^*B$. 
By shrinking $X_1$ if necessary, we can 
take a cyclic Galois covering $\pi_2: X_2 \to X_1$ of order $n_1$ 
which ramifies along the irreducible components of $B_1$ whose 
coefficients are equal to $1 - 1/n_1$.
Let $B_1'$ be the sum of all the other irreducible components.
Then we can show that the pair $(X_2, B_2)$ with $B_2 = \pi_2^*B_1'$
is again log terminal by a similar argument as in (1).

Since the set of integers arising from the coefficients of new pair 
$(X_2, B_2)$ is smaller, we 
obtain our assertion by induction.
\end{proof}

\begin{Expl}
(1) Let $X = \mathbb{C}^2/\mathbb{Z}_8(1,3)$, and
$f: Y \to X$ the minimal resolution.
Then $f^*K_X=K_Y+\frac 12 (C_1+C_2)$ for the exceptional divisors 
$C_1$ and $C_2$; 
$f$ is log crepant as a morphism from $(Y, \frac 12 (C_1+C_2))$ to 
$(X, 0)$.
Furthermore, let $g: Z \to Y$ be the blowing up at the point 
$y_0=C_1 \cap C_2$.
Then $g^*(K_Y+\frac 12 (C_1+C_2))=K_Z+\frac 12 (C'_1+C'_2)$, where 
$C'_i$ is the strict transform of $C_i$ for $i=1,2$.
Thus $f \circ g$ is log crepant as a morphism from 
$(Z, \frac 12 (C'_1+C'_2))$ to $(X, 0)$.

(2) Let $X=\mathbb{C}^2$, $D_i$ ($i=1,2$) the coordinate curves, 
and $f: Y \to X$ the weighted blowing up at 
the origin $x_0 = D_1 \cap D_2$ with weight $(1,2)$.
Then $Y$ has one ordinary double point, and 
$f^*(K_X+\frac 23(D_1+D_2))=K_Y+\frac 23(D'_1+D'_2)$, where 
$D'_i$ is the strict transform of $D_i$ for $i=1,2$.
Thus $f$ is log crepant as a morphism from $(Y, \frac 23(D'_1+D'_2))$ to 
$(X, \frac 23(D_1+D_2))$.
\end{Expl}

%%%%%%%%%%%%%%%%%%%%%%%%%%%%%%%%%%%%%%%%%%%%%%%%%%%%%%%%%%%%%%%%%%

\section{Recovery of varieties from categories}

The variety cannot be recovered from the derived category, but the canonical 
divisor can because the Serre functor is categorical:

\begin{Thm}\label{derived}
Let $(X,B)$ and $(Y,C)$ be pairs of projective varieties 
with $\mathbb{Q}$-divisors which satisfy the condition (*) in 
Definition~\ref{stack},
and let $\mathcal{X}$ and $\mathcal{Y}$ be the associated stacks.
Assume that there is an equivalence as triangulated categories:
\[
F: D^b(\text{Coh}(\mathcal{X})) \to D^b(\text{Coh}(\mathcal{Y})).
\]
Then the following hold:

(1) $\dim X = \dim Y$.

(2) $K_X+B$ (resp. $-(K_X+B)$) is nef if and only if $K_Y+C$ 
(resp. $-(K_Y+C)$) is nef.
Moreover, if this is the case, then the numerical Kodaira dimensions 
are equal
\[
\nu(X, \pm (K_X+B)) = \nu(Y, \pm (K_Y+C)).
\]

(3) If $K_X+B$ or $-(K_X+B)$ is big, then there are 
birational morphisms $\mu: W \to X$ and $\nu: W \to Y$ 
from a third projective variety $W$
such that $\mu^*(K_X+B) = \nu^*(K_Y+C)$.

(4) There is an isomorphism of big canonical rings as graded 
$\mathbb{C}$-algebras
\[
\bigoplus_{m \in \mathbb{Z}} H^0(X, \llcorner m(K_X+B) \lrcorner)
\to \bigoplus_{m \in \mathbb{Z}} H^0(Y, \llcorner m(K_Y+C) \lrcorner).
\]
\end{Thm}

\begin{proof}
By \cite{equi}, there exists an object on the product stack 
$e \in D^b(\text{Coh}(\mathcal{X} \times \mathcal{Y}))$ which is the 
kernel of the equivalence: $F$ is isomorphic to the integral functor 
$\Phi_{\mathcal{X} \to \mathcal{Y}}^e$ defined by 
$\Phi_{\mathcal{X} \to \mathcal{Y}}^e(a) = p_{2*}(p_1^*a \otimes e)$.
Then the proofs of the statements (1) through (3) are similar to those in 
\cite{DK}.

(4) This is due to Bridgeland and Caldararu.
Let $e' \in D^b(\text{Coh}(\mathcal{X} \times \mathcal{Y}))$ be the 
kernel of the quasi-inverse of $F$.
Let $\Delta_{\mathcal{X}}: \mathcal{X} \to \mathcal{X} \times \mathcal{X}$ and
$\Delta_{\mathcal{Y}}: \mathcal{Y} \to \mathcal{Y} \times \mathcal{Y}$ be
diagonal morphisms.
Then $e \boxtimes e' \in D^b(\text{Coh}(\mathcal{X} \times \mathcal{Y} \times 
\mathcal{X} \times \mathcal{Y}))$ is the kernel of an equivalence
\[
G: D^b(\text{Coh}(\mathcal{X} \times \mathcal{X})) 
\to D^b(\text{Coh}(\mathcal{Y} \times \mathcal{Y}))
\]
which satisfies
\[
G(\Delta_{{\mathcal{X}}*}\omega_{\mathcal{X}}^m) \cong 
\Delta_{{\mathcal{Y}}*}\omega_{\mathcal{Y}}^m
\]
for any integer $m$ by \cite{Caldararu}.
The conclusion follows from the following 
isomorphisms for $m, m' \in \mathbb{Z}$
\[
\begin{split}
&\text{Hom}_{\mathcal{X} \times \mathcal{X}}
(\Delta_{{\mathcal{X}}*}\omega_{\mathcal{X}}^m, 
\Delta_{{\mathcal{X}}*}\omega_{\mathcal{X}}^{m'})
\cong \text{Hom}_{\mathcal{X}}
(\Delta_{\mathcal{X}}^*\Delta_{{\mathcal{X}}*}\omega_{\mathcal{X}}^m, 
\omega_{\mathcal{X}}^{m'}) \\
&\cong \text{Hom}_{\mathcal{X}}
(\omega_{\mathcal{X}}^m \otimes 
\bigoplus_{p=0}^{\dim X} \Omega^p_{\mathcal{X}}[p], \omega_{\mathcal{X}}^{m'}) 
\cong \text{Hom}_{\mathcal{X}}
(\omega_{\mathcal{X}}^m, \omega_{\mathcal{X}}^{m'}).
\end{split}
\]
\end{proof}

The variety can be recovered from the category of sheaves:

\begin{Thm}\label{coh}
Let $(X,B)$ and $(Y,C)$ be pairs of projective varieties 
with $\mathbb{Q}$-divisors which satisfy the condition (*) in 
Definition~\ref{stack},
and let $\mathcal{X}$ and $\mathcal{Y}$ be the associated stacks.
Assume that there is an equivalence as abelian categories:
\[
F: \text{Coh}(\mathcal{X}) \to \text{Coh}(\mathcal{Y}).
\]
Then there exists an isomorphism $f: X \to Y$ such that $f_*B = C$.
Moreover, if $X$ is smooth and $B=0$, then there exists 
an invertible sheaf $L$ on $Y$ such that $F$ is isomorphic to the functor 
$f_*^L$ defined by $f_*^L(a) = f_*a \otimes L$.
\end{Thm}

\begin{proof}
The functor $F$ is extendable to an equivalence as triangulated categories
\[
\tilde F: D^b(\text{Coh}(\mathcal{X})) \to D^b(\text{Coh}(\mathcal{Y})).
\]
By \cite{equi}, there exists an object on the product stack 
$e \in D^b(\text{Coh}(\mathcal{X} \times \mathcal{Y}))$ such that 
$F$ is isomorphic to the integral functor 
$\Phi_{\mathcal{X} \to \mathcal{Y}}^e$.
The object $\mathcal{O}_x \in \text{Coh}(\mathcal{X})$ for $x \in X$
is simple, i.e, is non-zero and has no non-trivial subobject.
Then the image $p_{2*}(p_1^*\mathcal{O}_x \otimes e)$ is also simple and 
supported at a point.
Therefore, the support of $e$ is a graph of a bijection $f: X \to Y$, which is 
an isomorphism by the Zariski main theorem.
The coefficient of $B$ at any irreducible component is determined by the 
number of non-isomorphic simple objects above its generic point.
It follows that $C = f_*B$.
 
Under the additional assumption, we have 
\[
p_{2*}(p_1^*\mathcal{O}_x \otimes e) \cong \mathcal{O}_{f(x)}
\]
for any point $x \in X$.
Therefore, $e$ is an invertible sheaf on the graph of $f$.
\end{proof}

%%%%%%%%%%%%%%%%%%%%%%%%%%%%%%%%%%%%%%%%%%%%%%%%%%%%%%%%%%%%%%%%%

\section{Toroidal case}

We confirm our conjecture in the case of toroidal varieties.

\begin{Defn}
Let $X$ be a normal variety and $\bar B$ a reduced divisor.
The pair $(X, \bar B)$ is said to be {\em toroidal} if, for each $x \in X$, 
there exists a toric variety $(P_x, Q_x)$ with a point $t$, 
called a {\em toric chart}, such that
the completion $(\hat X_x, \hat{\bar B}_x)$ at $x$ is isomorphic to the 
completion $(\hat{(P_x)}_t, \hat{(Q_x)}_t)$ at $t$.  
A toroidal pair is said to be {\em quasi-smooth} if it has only quotient 
singularities.  Such a pair always satisfies the condition (*) in 
Definition~\ref{stack}.

A morphism $f: (X, \bar B) \to (Y, \bar C)$ of toroidal pairs is said to be 
{\em toroidal} if, for any point $x \in X$ and any toric chart 
$(P'_y, Q'_y; t')$ at $y = f(x) \in Y$, 
there exists a toric chart $(P_x, Q_x; t)$ at $x$ and a 
toric morphism $g_x:(P_x, Q_x) \to (P'_y, Q'_y)$ such that $g_x(t)=t'$ and 
the completion $\hat f_x: (\hat X_x, \hat{\bar B}_x) \to 
(\hat Y_y, \hat{\bar C}_y)$ is isomorphic to the 
completion $\hat{(g_x)}_t: (\hat{(P_x)}_t, \hat{(Q_x)}_t)
\to (\hat{(P'_y)}_{t'}, \hat{(Q'_y)}_{t'})$  
under the toric chart isomorphisms.

A toroidal morphism $f: (X, B) \to (Y, C)$ is said to be a {\em divisorial
contraction} if it is a projective birational morphism and 
such that the exceptional locus is a prime divisor.
A pair of toroidal morphisms $\phi: (X, B) \to (Z, D)$ and $\psi: 
(Y, C) \to (Z, D)$ is said to be a {\em flip} if they are projective 
birational morphisms such that the codimensions of their exceptional loci are 
at least $2$, the composite map $(\psi)^{-1} \circ \phi$ is not an 
isomorphism, and that the relative Picard numbers $\rho(X/Z)$ and $\rho(Y/Z)$ 
are equal to $1$.
\end{Defn}

\begin{Thm}\label{toroidal}
Let $X$ and $Y$ be quasi-smooth toroidal varieties, 
$B$ and $C$ be effective toroidal $\mathbb{Q}$-divisors
on $X$ and $Y$, respectively, whose coefficients are contained in the set 
$\{1 - 1/r \vert r \in \mathbb{Z}_{>0}\}$, and let $\mathcal{X}$ and 
$\mathcal{Y}$ be smooth Deligne-Mumford stacks associated to the 
pairs $(X, B)$ and $(Y, C)$, respectively. 
Assume that one of the following holds.

(1) $X = Y$ and $K_X + B \ge K_Y + C$.

(2) There is a toroidal divisorial contraction $f: X \to Y$
such that $C = f_*B$ and 
\[
K_X+B \ge f^*(K_Y + C).
\]

(3) There is a toroidal flip $f: X \to Z \leftarrow Y$ 
such that $C = f_*B$ and 
\[
\mu^*(K_X+B) \ge \nu^*(K_Y + C)
\]
for common toroidal resolutions $\mu: W \to X$ and $\nu: W \to Y$ with
$f = \nu \circ \mu^{-1}$.

(4) There is a toroidal divisorial contraction $f: Y \to X$ 
such that $B = f_*C$ and 
\[
f^*(K_X+B) \ge K_Y + C.
\]

Setting $Z = Y$ in the cases (1) and (2) and $Z = X$ in the case (4), let
$\mathcal{W} = (\mathcal{X} \times_Z \mathcal{Y})\tilde{}$ be the 
normalization of the fiber product with the natural morphisms
$\tilde{\mu}: \mathcal{W} \to \mathcal{X}$ and 
$\tilde{\nu}: \mathcal{W} \to \mathcal{Y}$.
Then the functor
\[
F = \tilde{\mu}_*\tilde{\nu}^*: 
D^b(\text{Coh}(\mathcal{Y})) \to D^b(\text{Coh}(\mathcal{X})).
\] 
is fully faithful. 
Moreover, if the inequality between the log canonical divisors become 
an equality, then $F$ is an equivalence of triangulated categories.
\end{Thm}

\begin{proof}
Since the set of all the point sheaves span the derived categories, 
our assertion can be proved analytic locally.
Hence we may assume that $X, Y, Z$ are toric varieties, $Z$ is affine, 
and that $f$ is a toric map.
We note that the toroidal structure is only necessary 
near the exceptional locus of the given birational maps for the application
of the theorem.

It is sufficient to prove the first assertion.
Indeed, if the pull backs of $K_X+B$ and $K_Y + C$ are equal, 
then $F$ is automatically an equivalence, 
because it commutes with the Serre functors (\cite{Bridgeland}).

Let $N$ be a lattice in $\mathbb{R}^n$ which contains $\mathbb{Z}^n$,
$M$ its dual lattice.
Assume that the standard vectors 
$v_1=(1,0, \dots, 0)$, \dots, $v_n=(0,\dots,0,1)$ are primitive in $N$.
Let $\{v_i^*\}_{1 \le i \le n}$ be the dual basis of $\{v_i\}_{1 \le i \le n}$.

\vskip 1pc

(1) Let $X = Y = Z$ be the affine toric variety corresponding to
the cone $\sigma \subset N_{\mathbb{R}}$ generated by 
$v_1, \dots, v_n$.
Let $D_i$ be the prime divisors on $X$ corresponding to the $v_i$.
Let $r_i$ and $s_i$ be positive integers attached to the divisors
$D_i$ for $1 \le i \le n$ such that we have 
$B = \sum_i (1 - \frac 1{r_i})D_i$ and 
$C = \sum_i (1 - \frac 1{s_i})D_i$.
By assumption, we have $r_i \ge s_i$ for all $i$.

The stacks $\mathcal{X}$ and $\mathcal{Y}$ are defined by the 
coverings $\pi_X: \tilde X \to X$ 
and $\pi_Y: \tilde Y \to Y$ which correspond to the sublattices
$N_X$ and $N_Y$ of $N$ generated by the $r_iv_i$ and $s_iv_i$ 
for $1 \le i \le n$, respectively.
The stack $\mathcal{W}$ is defined 
by the covering $\pi_W: \tilde W \to X$ which correspond to the sublattice
$N_W = N_X \cap N_Y \subset N$.
Let $M_X$, $M_Y$ and $M_W$ be the dual lattices of $N_X$, $N_Y$ and $N_W$, 
respectively.
The morphisms $\tilde{\mu}: \mathcal{W} \to \mathcal{X}$ and 
$\tilde{\nu}: \mathcal{W} \to \mathcal{Y}$ are given by 
inclusions $M_X \subset M_W$ and $M_Y \subset M_W$, respectively.
$M_X$ and $M_Y$ are generated by the $\frac 1{r_i}v_i^*$ and 
$\frac 1{s_i}v_i^*$, respectively.
The dual lattice $M_W$ is equal to 
the sum $M_X+M_Y \subset M_{\mathbb{Q}}$ generated by the $\frac 1{t_i}v_i^*$
with $t_i = \text{LCM}(r_i,s_i)$.

Let $m = \sum_{i=1}^n \frac{m_i}{s_i}v_i^* \in M_Y$ 
for $m_i \in \mathbb{Z}$ be a monomial 
which generates an invertible sheaf 
$L = \mathcal{O}_{\mathcal{Y}}(- \sum_{i=1}^n \frac{m_i}{s_i}D_i)$ 
on $\mathcal{Y}$.
Then the invertible sheaf $\tilde{\nu}^*L$ on $\mathcal{W}$ is 
generated by the same monomial $m$, hence $F(L)$ on $\mathcal{X}$
is by the generator
\[
F(m) = \sum_{i=1}^n \frac{\ulcorner \frac{m_ir_i}{s_i} \urcorner}{r_i}v_i^*
\]
of $(m + \check{\sigma}) \cap M_X$. 
Let $L'$ be another invertible sheaf with 
the corresponding monomials
\[
m'=\sum_{i=1}^n \frac{m'_i}{s_i}v_i^*, \quad  
F(m') = \sum_{i=1}^n \frac{\ulcorner \frac{m'_ir_i}{s_i} \urcorner}
{r_i}v_i^*. 
\]
We have 
\[
\begin{split}
&\text{Hom}(L', L) = \mathbb{C}^{(m - m' + \check{\sigma}) \cap M} \\
&\text{Hom}(F(L'), F(L)) = \mathbb{C}^{(F(m) - F(m') + \check{\sigma}) 
\cap M}. 
\end{split}
\]
We calculate 
\[
0 \le \ulcorner \frac{m_i-m'_i}{s_i} \urcorner r_i 
+ \ulcorner \frac{m'_ir_i}{s_i} \urcorner - \frac{m_ir_i}{s_i}
\le (1 - \frac 1{s_i})r_i + (1 - \frac 1{s_i}) < r_i 
\]
because $s_i \le r_i$.
Thus
\[
0 \le \ulcorner \frac{m_i-m'_i}{s_i} \urcorner r_i 
+ \ulcorner \frac{m'_ir_i}{s_i} \urcorner
- \ulcorner \frac{m_ir_i}{s_i} \urcorner < r_i
\]
hence
\[
\ulcorner \frac{m_i-m'_i}{s_i} \urcorner
= \ulcorner \frac{\ulcorner \frac{m_ir_i}{s_i} \urcorner 
- \ulcorner \frac{m'_ir_i}{s_i} \urcorner}{r_i} \urcorner.
\]
Therefore, the homomorphism
\[
\text{Hom}(L', L) \to \text{Hom}(F(L'), F(L))
\]
is bijective.

%%%%%%%%%%%%%%%%%%%%%%%%%%%%%%%%%%%%%%%%%%%%%%%%%%%%%%%%%%%%%%%%%

\vskip 1pc

(2) Let $Y = Z$ be the affine toric variety corresponding to
the cone $\sigma \subset N_{\mathbb{R}}$ generated by 
$v_1, \dots, v_n$ as in (1).
Let $v_{n+1}=(a_1,\dots,a_n) \in N$ be a primitive vector 
such that $a_i > 0$ for $1 \le i \le n'$ and $a_i=0$ for 
$n' < i \le n$.
We assume that $n' \ge 2$.
If we put $a_{n+1}=-1$, then we have
\[
\sum_{i=1}^{n+1} a_iv_i=0.
\]
Let $\sigma_{i_0}$ ($1 \le i_0 \le n'$) be the subcones of $\sigma$
generated by the $v_i$ for $1 \le i \le n$ with $i \ne i_0$.
Then there is a decomposition $\sigma = \bigcup_{i=1}^{n'} \sigma_i$.
Let $f: X \to Y$ be the corresponding projective birational toric morphism,
let $D_i$ be the prime divisors on $X$ corresponding to the $v_i$ for 
$1 \le i \le n+1$, and $D'_i$ their strict transforms on $Y$. 
$D_{n+1}$ is the exceptional divisor of $f$, and we put $D'_{n+1}=0$.
Then we have $f^*D'_i = D_i + a_iD_{n+1}$, hence
$D_i \equiv -a_iD_{n+1}$.
The divisors $D_i$ for $1 \le i \le n'$ are positive for $f$.

Let $r_i$ be positive integers attached to the divisors
$D_i$ for $1 \le i \le n+1$ such that we have 
$B = \sum_i (1 - \frac 1{r_i})D_i$.
Our condition 
\[
K_X+\sum_{i=1}^{n+1} \frac{r_i-1}{r_i}D_i
\ge f^*(K_Y+\sum_{i=1}^n \frac{r_i-1}{r_i}D'_i).
\]
is equivalent to 
\[
\sum_{i=1}^{n+1} \frac{a_i}{r_i} \ge 0
\]
because $K_X+\sum_{i=1}^{n+1} D_i = f^*(K_Y+\sum_{i=1}^n D'_i)$.

The toric variety $X$ is covered by the affine toric varieties
$X_{\sigma_{i_0}}$ for $1 \le i_0 \le n'$. 
The stacks $\mathcal{X}$ and $\mathcal{Y}$ are defined by the 
coverings $\pi_{X,i_0}: \tilde X_{\sigma_{i_0}} \to X_{\sigma_{i_0}}$ 
and $\pi_Y: \tilde Y \to Y$ 
which correspond to the sublattices
$N_{i_0}$ of $N$ generated by the $r_iv_i$ for $1 \le i \le n+1$
with $i \ne i_0$ and $N_{n+1}$ generated by the $r_iv_i$ for $1 \le i \le n$,
respectively.
The stack $\mathcal{W}$ above $X$ is defined 
by the coverings $\pi_{W,i_0}: \tilde W_{\sigma_{i_0}} \to 
X_{\sigma_{i_0}}$ which correspond to the sublattices
$N_{i_0, n+1} = N_{i_0} \cap N_{n+1} \subset N$.
We note that the morphism of stacks
$\tilde{\mu}: \mathcal{W} \to \mathcal{X}$ is not necessarily an isomorphism.
If $1 \le i_0 \le n$, then the dual lattice $M_{i_0}$ of $N_{i_0}$  
is generated by $\frac 1{r_i}v_i^* - \frac{a_i}{a_{i_0}r_i}v_{i_0}^*$
for $1 \le i \le n$ with $i \ne i_0$ and $\frac 1{a_{i_0}r_{n+1}}v_{i_0}^*$.
The dual lattice $M_{n+1}$ of $N_{n+1}$ is generated by 
$\frac 1{r_i}v_i^*$ for $1 \le i \le n$.
The dual lattice $M_{i_0, n+1}$ of $N_{i_0, n+1}$ is equal to 
the sum $M_{i_0}+M_{n+1} \subset M_{\mathbb{Q}}$.

Let $m  = \sum_{i=1}^n \frac{k_i}{r_i}v_i^* \in M_{n+1}$
with $k_i \in \mathbb{Z}$ be a monomial which corresponds to an 
invertible sheaf $L$ on $\mathcal{Y}$.
Then the invertible sheaf $F(L)$ on $\mathcal{X}$
is given by the generators $m_{i_0}$ of 
$(m + \check{\sigma}_{i_0}) \cap M_{i_0}$ given by
\[
m_{i_0} = \sum_{1 \le i \le n, i\ne i_0} \frac{k_i}{r_i}v_i^* +
(\frac{k_{n+1}}{a_{i_0}r_{n+1}}-\sum_{1 \le i \le n, i\ne i_0}
\frac{a_ik_i}{a_{i_0}r_i})v_{i_0}^*
\]
where $k_{n+1}$ is a smallest integer such that
\[
\frac{k_{n+1}}{a_{i_0}r_{n+1}}-\sum_{1 \le i \le n, i\ne i_0}
\frac{a_ik_i}{a_{i_0}r_i} \ge \frac{k_{i_0}}{r_{i_0}}
\]
that is
\[
k_{n+1}=\ulcorner r_{n+1} \sum_{i=1}^n \frac{a_ik_i}{r_i} \urcorner.
\]
We note that $k_{n+1}$ is independent of $i_0$.
Let $L'$ be another invertible sheaf corresponding to the 
monomial $m'=\sum_{i=1}^n \frac{k'_i}{r_i}v_i^*$, and let
\[
m'_{i_0} = \sum_{1 \le i \le n, i\ne i_0} \frac{k'_i}{r_i}v_i^*
-\sum_{1 \le i \le n+1, i\ne i_0}\frac{a_ik'_i}{a_{i_0}r_i}v_{i_0}^*, \quad
k_{n+1}'=\ulcorner r_{n+1} \sum_{i=1}^n \frac{a_ik'_i}{r_i} \urcorner
\]
where we note that $a_{n+1}=-1$.

Since $f$ is birational and $\mathcal{H}om(L', L)$ is torsion free, 
the natural homomorphism
\[
\text{Hom}(L', L) \to \text{Hom}(F(L'), F(L)) 
\]
is injective.  
It is also surjective because the complement of the indeterminacy locus of 
$f^{-1}$ in $Y$ has codimension at least $2$.
Since $\text{Hom}^p(L', L)=0$ for $p>0$, it is sufficient to prove that
$\text{Hom}^p(F(L'), F(L))=0$ for $p>0$ in order to prove that 
$F$ is fully faithful.

The invertible sheaf $\mathcal{H}om(F(L'),F(L))$ is given by the monomials
\[
m_{i_0}-m'_{i_0} = \sum_{1 \le i \le n, i\ne i_0} \frac{k_i-k'_i}{r_i}v_i^* 
-\sum_{1 \le i \le n+1, i\ne i_0} \frac {a_i(k_i-k'_i)}{a_{i_0}r_i}v_{i_0}^*.
\]
We consider the divisorial reflexive sheaf 
\[
\mathcal{H}om(F(L'),F(L))_X = 
\mathcal{O}_X(- \sum_{i=1}^{n+1} \ulcorner \frac{k_i-k'_i}{r_i} \urcorner D_i)
\]
on $X$.  
Since $f^*D'_i \equiv D_i + a_iD_{n+1}$ for $1 \le i \le n$,
we have 
\[
- \sum_{i=1}^{n+1} \ulcorner \frac{k_i-k'_i}{r_i} \urcorner D_i
\equiv \sum_{i=1}^{n+1} a_i \ulcorner \frac{k_i-k'_i}{r_i} \urcorner D_{n+1}.
\]
We calculate 
\[
\begin{split}
&\sum_{i=1}^{n+1} \frac{a_i(k_i-k'_i)}{r_i} \\
&= \sum_{i=1}^n \frac{a_i(k_i-k'_i)}{r_i} 
- \frac{\ulcorner r_{n+1} \sum_{i=1}^n \frac{a_ik_i}{r_i} \urcorner}{r_{n+1}}
+ \frac{\ulcorner r_{n+1} \sum_{i=1}^n \frac{a_ik'_i}{r_i} \urcorner}{r_{n+1}} 
\\
&< \frac 1{r_{n+1}} \le \sum_{i=1}^n \frac{a_i}{r_i}.
\end{split}
\]
Therefore we have 
\[
\sum_{i=1}^{n+1} a_i \ulcorner \frac{k_i-k'_i}{r_i} \urcorner 
\le \sum_{i=1}^{n+1} \frac{a_i (k_i-k'_i)}{r_i} + 
\sum_{i=1}^n \frac{a_i(r_i-1)}{r_i}
< \sum_{i=1}^n a_i.
\]
On the other hand, we have 
$K_X + D_{n+1} \equiv (\sum_{i=1}^n a_i)D_{n+1}$.
Since $- D_{n+1}$ is an $f$-ample $\mathbb{Q}$-divisor, 
we conclude by the vanishing theorem (\cite{KMM}~Theorem~1.2.5) that 
\[
\text{Hom}^p(F(L'),F(L)) \cong H^p(\mathcal{H}om(F(L'),F(L))_X)=0
\]
for $p > 0$.

%%%%%%%%%%%%%%%%%%%%%%%%%%%%%%%%%%%%%%%%%%%%%%%%%%%%%%%%%%%%%%%%%%%

\vskip 1pc

(3) Let $v_{n+1}=(a_1,\dots,a_n) \in N$ be a primitive vector 
such that $a_i > 0$ for $1 \le i \le n'$, $a_i=0$ for 
$n' < i \le n''$ and $a_i < 0$ for $n'' < i \le n$.
We assume that $2 \le n'$ and $n'' < n$.
If we put $a_{n+1}=-1$, then we have
\[
\sum_{i=1}^{n+1} a_iv_i=0.
\]
Let $\langle v_1, \dots, v_{n+1} \rangle$ be the convex cone  
generated by the $v_i$ for $1 \le i \le n+1$,
and let $\sigma_{i_0}$ ($1 \le i_0 \le n +1$) be the subcones 
generated by the $v_i$ for $1 \le i \le n+1$ with $i \ne i_0$.
Then there are two decompositions:
\[
\langle v_1, \dots, v_{n+1} \rangle
= \bigcup_{i=1}^{n'} \sigma_i 
= \bigcup_{i=n''+1}^{n+1} \sigma_i. 
\]

Let $Z$ be an affine toric variety corresponding to the lattive $N$ and 
the cone $\langle v_1, \dots, v_{n+1}\rangle$.
Let $g: X \to Z$ and $h: Y \to Z$ be projective birational toric morphisms
corresponding to these subdivisions, and let
$f = h \circ g^{-1}: X -\to Y$ be the composite proper birational map.
Let $D_i$ be the prime divisors on $X$ corresponding to the $v_i$ for 
$1 \le i \le n+1$, and $D'_i$ their strict transforms on $Y$. 
Then we have
$D_i \equiv -a_iD_{n+1}, \quad D'_i \equiv -a_iD'_{n+1}$.
The divisors $D_i$ ($1 \le i \le n'$) are positive for $g$, 
and the $D_i$ ($n'' < i \le n+1$) negative.

Let $v_{n+2}=\sum_{i=1}^{n'}\lambda a_iv_i=
\sum_{i=n''+1}^{n+1}\lambda (-a_i)v_i$, where $\lambda$ is the smallest 
positive number such that $v_{n+2}$ becomes a primitive vector in $N$.
We define the subcones
$\sigma_{i_0i_1}$ for $1 \le i_0 \le n'$ and $n'' < i_1 \le n+1$ to be 
the ones generated by the $v_i$ for $1 \le i \le n+1$ with $i \ne i_0,i_1$
and $v_{n+2}$.
Then we have two decompositions
\[
\sigma_{i_1} = \bigcup_{i=1}^{n'} \sigma_{ii_1}, \quad
\sigma_{i_0} = \bigcup_{i=n''+1}^{n+1} \sigma_{i_0i}
\]
so that 
\[
\langle v_1, \dots, v_{n+1} \rangle
= \bigcup_{1 \le i \le n', n''+1 \le j \le n+1} \sigma_{ij}. 
\]
Let $W$ be the toric variety corresponding to this subdivision with
natural projective birational morphisms $\mu: W \to X$ and $\nu: W \to Y$ :
\[
\begin{CD}
W @>{\nu}>> Y \\
@V{\mu}VV @VVhV \\
X @>>g> Z.
\end{CD}
\] 
Let $D''_i$ be the prime divisors on $W$ corresponding to the $v_i$ for 
$1 \le i \le n+2$. 
We have
\[
\begin{split}
&\mu^*D_i = \begin{cases} D''_i &\text{ if } 1 \le i \le n' \\
D''_i + \lambda (- a_i)D''_{n+2} &\text{ if } n'' < i \le n+1 
\end{cases} \\
&\nu^*D'_i = \begin{cases} D''_i + \lambda a_iD''_{n+2} &\text{ if } 
1 \le i \le n' \\
D''_i &\text{ if } n'' < i \le n+1. 
\end{cases}
\end{split}
\]
If $r_1, \dots, r_{n+1}$ are the positive integers attached to the divisors
$D_i$, then our condition 
\[
\mu^*(K_X+\sum_{i=1}^{n+1} \frac{r_i-1}{r_i}D_i)
\ge \nu^*(K_Y+\sum_{i=1}^{n+1} \frac{r_i-1}{r_i}D'_i).
\]
is equivalent to 
\[
\sum_{i=1}^{n+1} \frac{a_i}{r_i} \ge 0
\]
because $\mu^*(K_X+\sum_{i=1}^{n+1} D_i)
= \nu^*(K_Y+\sum_{i=1}^{n+1} D'_i)
= K_W + \sum_{i=1}^{n+2} D''_i$.

The toric varieties $X$ and $Y$ are covered by the affine toric varieties
$X_{\sigma_{i_0}}$ for $1 \le i_0 \le n'$ and $Y_{\sigma_{i_0}}$ for 
$n'' < i_0 \le n+1$, respectively. 
The stacks $\mathcal{X}$ and $\mathcal{Y}$ are defined by the 
coverings $\pi_{X,i_0}: \tilde X_{\sigma_{i_0}} \to X_{\sigma_{i_0}}$ 
and $\pi_{Y,i_0}: \tilde Y_{\sigma_{i_0}} \to Y_{\sigma_{i_0}}$
which correspond to the sublattices
$N_{i_0}$ of $N$ generated by the $r_iv_i$ for $1 \le i \le n+1$
with $i \ne i_0$.
$W$ is covered by $W_{\sigma_{i_0i_1}}$ for $1 \le i_0 \le n'$ and 
$n'' < i_1 \le n+1$, and the stack $\mathcal{W}$ is defined 
by the coverings $\pi_{W,i_0i_1}: \tilde W_{\sigma_{i_0i_1}} \to 
W_{\sigma_{i_0i_1}}$ which correspond to the sublattices
$N_{i_0i_1} = N_{i_0} \cap N_{i_1} \subset N$.
The morphisms $\mu: W \to X$ and $\nu: W \to Y$ are covered by
morphisms $\tilde{\mu}: \mathcal{W} \to \mathcal{X}$ and 
$\tilde{\nu}: \mathcal{W} \to \mathcal{Y}$ of stacks.
If $1 \le i_0 \le n$, then the dual lattice $M_{i_0}$ of $N_{i_0}$  
is generated by $\frac 1{r_i}v_i^* - \frac{a_i}{a_{i_0}r_i}v_{i_0}^*$
for $1 \le i \le n$ with $i \ne i_0$ and $\frac 1{a_{i_0}r_{n+1}}v_{i_0}^*$.
The dual lattice $M_{n+1}$ of $N_{n+1}$ is generated by 
$\frac 1{r_i}v_i^*$ for $1 \le i \le n$.
The dual lattice $M_{i_0i_1}$ of $N_{i_0i_1}$ is equal to 
the sum $M_{i_0}+M_{i_1} \subset M_{\mathbb{Q}}$.

Let $L = \mathcal{O}_{\mathcal{Y}}(\sum_{i=1}^{n+1} \frac{k_i}{r_i}D'_i)$
be an invertible sheaf on $\mathcal{Y}$, where $k_i$ are 
integers.
We have 
\[
\begin{split}
&\tilde{\nu}^*L = 
\mathcal{O}_{\mathcal{W}}(\sum_{i=1}^{n+1} \frac{k_i}{r_i}D''_i
+\sum_{i=1}^{n'} \frac{\lambda a_ik_i}{r_i}D''_{n+2}) \\
&= \tilde{\mu}^*\mathcal{O}_{\mathcal{X}}(\sum_{i=1}^{n+1} 
\frac{k_i}{r_i}D_i) \otimes \mathcal{O}_{\mathcal{W}}
(\sum_{i=1}^{n+1} \frac{\lambda a_ik_i}{r_i}D''_{n+2}).
\end{split}
\]
Since  
\[
K_W+\sum_{i=1}^{n+1} \frac{r_i-1}{r_i}D''_i+D''_{n+2}
=\mu^*(K_X+\sum_{i=1}^{n+1} \frac{r_i-1}{r_i}D_i)
-\sum_{i=n''+1}^{n+1}\frac{\lambda a_i}{r_i}D''_{n+2}
\]
we have $R^p\tilde{\mu}_*\tilde{\nu}^*L=0$ for $p > 0$,
if $\sum_{i=1}^{n+1} \frac{a_ik_i}{r_i} < 
- \sum_{i=n''+1}^{n+1}\frac{a_i}{r_i}$ by \cite{KMM}~Theorem~1.2.5.
Thus if 
\begin{equation}\label{range}
0 \le \sum_{i=1}^{n+1} \frac{a_ik_i}{r_i} 
< - \sum_{i=n''+1}^{n+1}\frac{a_i}{r_i}
\end{equation}
then we have $F(L) = \mathcal{O}_{\mathcal{X}}(\sum_{i=1}^{n+1} 
\frac{k_i}{r_i}D_i)$.

Let $L' = \mathcal{O}_{\mathcal{Y}}(\sum_{i=1}^{n+1} \frac{k'_i}{r_i}D'_i)$
be another invertible sheaf. 
We assume that both $L$ and $L'$ are in the range (\ref{range}).
Since $f$ is isomorphic in codimension $1$, the natural homomorphism 
\[
\text{Hom}(L', L) \to \text{Hom}(F(L'), F(L))
\]
is bijective.
We shall prove that their higher Hom's vanish on both $X$ and $Y$.

We have $\mathcal{H}om(L', L)_Y 
= \mathcal{O}_Y(\sum_{i=1}^{n+1} \llcorner \frac{k_i-k'_i}{r_i} 
\lrcorner D'_i)$ and
\[
\sum_{i=1}^{n+1} \llcorner \frac{k_i-k'_i}{r_i} \lrcorner D'_i
\equiv - (\sum_{i=1}^{n+1} a_i \llcorner \frac{k_i-k'_i}{r_i} \lrcorner)
D'_{n+1}.
\]
Since $- \sum_{i=1}^{n+1} a_i\frac{k_i-k'_i}{r_i}
> \sum_{i=n''+1}^{n+1} \frac{a_i}{r_i}$, we have 
\[
- \sum_{i=1}^{n+1} a_i \llcorner \frac{(k_i-k'_i)}{r_i} \lrcorner
> \sum_{i=n''+1}^{n+1} \frac{a_i}{r_i} + 
\sum_{i=n''+1}^{n+1} a_i\frac{r_i-1}{r_i}
= \sum_{i=n''+1}^{n+1} a_i.
\]
Since $K_Y + \sum_{i=1}^{n'} D'_i \equiv - \sum_{i=n''+1}^{n+1}D'_i
\equiv (\sum_{i=n''+1}^{n+1} a_i)D'_{n+1}$ and $D'_{n+1}$ is ample,
we conclude that
\[
\text{Hom}^p(L', L)=H^p(Y, \mathcal{H}om(L', L)_Y) = 0
\]
for $p>0$ by \cite{KMM}~Theorem~1.2.5.

Similarly, we have $\mathcal{H}om(F(L'), F(L))_X 
= \mathcal{O}_X(\sum_{i=1}^{n+1} \llcorner \frac{k_i-k'_i}{r_i} \lrcorner 
D_i)$ and 
\[
\sum_{i=1}^{n+1} \llcorner \frac{k_i-k'_i}{r_i} \lrcorner D_i
\equiv - (\sum_{i=1}^{n+1} a_i \llcorner \frac{(k_i-k'_i)}{r_i} \lrcorner) 
D_{n+1}.
\]
Since $\sum_{i=1}^{n+1} \frac{a_i(k_i-k'_i)}{r_i}
> \sum_{i=n''+1}^{n+1} \frac{a_i}{r_i}$, we have 
\[
\begin{split}
&\sum_{i=1}^{n+1} a_i \llcorner \frac{k_i-k'_i}{r_i} \lrcorner
> \sum_{i=n''+1}^{n+1} \frac{a_i}{r_i} - \sum_{i=1}^{n'} a_i\frac{r_i-1}{r_i} 
\\
&\ge - \sum_{i=1}^{n'} \frac{a_i}{r_i} - \sum_{i=1}^{n'} a_i\frac{r_i-1}{r_i}
= - \sum_{i=1}^{n'} a_i.
\end{split}
\]
Since $K_X + \sum_{i=n''+1}^{n+1} D_i \equiv - \sum_{i=1}^{n'}D_i
\equiv (\sum_{i=1}^{n'} a_i)D_{n+1}$ and $-D_{n+1}$ is ample,
we conclude that
\[
\text{Hom}^p(F(L'), F(L))=H^p(X, \mathcal{H}om(F(L'), F(L))_X) = 0
\]
for $p>0$ by \cite{KMM}~Theorem~1.2.5.
Therefore, the proof of the fully faithfulness of $F$ 
is reduced to the following Lemma~\ref{span} by using \cite{Bridgeland}.

%%%%%%%%%%%%%%%%%%%%%%%%%%%%%%%%%%%%%%%%%%%%%%%%%%%%%%%%%%%%%%%%%

\vskip 1pc

(4) We interchange $X$ and $Y$, and use the notation of (2).
The inequality $K_X+B \le f^*(K_Y+C)$ implies that 
\[
\sum_{i=1}^n \frac{a_i}{r_i} \le \frac 1{r_n}.
\]
Let $L = \mathcal{O}_{\mathcal{X}}(\sum_{i=1}^{n+1} \frac{k_i}{r_i}D_i)$
be an invertible sheaf on $\mathcal{X}$, where $k_i$ are 
integers.
Since $f^*D'_i=D_i + a_iD_{n+1}$, we have
\[
\tilde{\mu}^*L = 
\mathcal{O}_{\mathcal{W}}(\sum_{i=1}^{n+1} \frac{k_i}{r_i}D''_i) 
= \tilde{\nu}^*\mathcal{O}_{\mathcal{Y}}(\sum_{i=1}^n \frac{k_i}{r_i}D'_i) 
\otimes \mathcal{O}_{\mathcal{W}}
(-\sum_{i=1}^{n+1} \frac{a_ik_i}{r_i}D''_{n+1})
\]
where we note that $a_{n+1}=-1$.
Since  
\[
K_X+\sum_{i=1}^n \frac{r_i-1}{r_i}D_i + D_{n+1}
=f^*(K_Y+\sum_{i=1}^n \frac{r_i-1}{r_i}D'_i)
+\sum_{i=1}^n\frac{a_i}{r_i}D_{n+1}
\]
we have $R^p\tilde{\nu}_*\tilde{\mu}^*L=0$ for $p > 0$ if 
$- \sum_{i=1}^{n+1} \frac{a_ik_i}{r_i} < \sum_{i=1}^n \frac{a_i}{r_i}$.
Thus if 
\begin{equation}\label{range2}
0 \le - \sum_{i=1}^{n+1} \frac{a_ik_i}{r_i} 
< \sum_{i=1}^n\frac{a_i}{r_i}
\end{equation}
then we have $F(L) = \mathcal{O}_{\mathcal{Y}}(\sum_{i=1}^n 
\frac{k_i}{r_i}D'_i)$.

Let $L' = \mathcal{O}_{\mathcal{X}}(\sum_{i=1}^{n+1} \frac{k'_i}{r_i}D_i)$
be another invertible sheaf, so that we have $\mathcal{H}om(L', L)_X 
= \mathcal{O}_X(\sum_{i=1}^{n+1} \llcorner \frac{k_i-k'_i}{r_i} 
\lrcorner D_i)$.
We assume that both $L$ and $L'$ are in the range (\ref{range2}).
Since we have $\sum_{i=1}^n \frac{a_i}{r_i} \le \frac 1{r_n}$, 
the homomorphism
\[
\text{Hom}(L', L) \to \text{Hom}(F(L'), F(L))
\]
is bijective. 

We have
\[
\sum_{i=1}^{n+1} \llcorner \frac{k_i-k'_i}{r_i} \lrcorner D_i
\equiv - (\sum_{i=1}^{n+1} a_i \llcorner \frac{k_i-k'_i}{r_i} \lrcorner) 
D_{n+1}.
\]
Since $- \sum_{i=1}^{n+1} a_i\frac{k_i-k'_i}{r_i}
< \sum_{i=1}^n \frac{a_i}{r_i}$, we have 
\[
- \sum_{i=1}^{n+1} a_i \llcorner \frac{(k_i-k'_i)}{r_i} \lrcorner
< \sum_{i=1}^n \frac{a_i}{r_i} + \sum_{i=1}^n a_i\frac{r_i-1}{r_i}
= \sum_{i=1}^n a_i.
\]
Since $K_X + D_{n+1} \equiv - \sum_{i=1}^n D_i
\equiv (\sum_{i=1}^n a_i)D_{n+1}$
and $D_{n+1}$ is ample, we conclude that
\[
\text{Hom}^p(L', L)=H^p(X, \mathcal{H}om(L', L)_X) = 0
\]
for $p>0$.
Therefore, $F$ is fully faithful by the following Lemma~\ref{span'}. 
\end{proof}

\begin{Lem}\label{span}
The set of all the invertible sheaves on $\mathcal{Y}$
in the range \ref{range} is a spanning class of the category 
$D^b(\text{Coh}(\mathcal{Y}))$.
\end{Lem}

\begin{proof}
Let $\Omega$ be the set of all such invertible sheaves.
We shall prove that the full subcategory $\mathcal{D}$ spanned by $\Omega$
coincides with $D^b(\text{Coh}(\mathcal{Y}))$.
Let $i_p$ ($1 \le p \le t$) be integers such that 
$n'' < i_1 < \dots < i_t \le n+1$ and $t < n-n''$.
Let $i_0$ be an integer such that $n'' < i_0 \le n+1$ and 
$i_0 \ne i_p$ for any $p$.
Let $S$ be the intersection of the Cartier divisors $\frac 1{r_{i_p}}D'_{i_p}$ 
for $1 \le p \le t$ on $\mathcal{Y}$.
We call $S$ a stratum.
We have an exact Koszul complex
\[
0 \to \mathcal{O}_{\mathcal{Y}}(- \sum_{p=1}^t \frac 1{r_{i_p}}D'_{i_p})
\to \dots \to 
\bigoplus_{p=1}^t \mathcal{O}_{\mathcal{Y}}(- \frac 1{r_{i_p}}D'_{i_p})
\to \mathcal{O}_{\mathcal{Y}} \to \mathcal{O}_S \to 0.
\]
Let $k_i$ ($1 \le i \le n+1$, $i \ne i_0$) be any choice of integers, and let
$k_{i_0}$ be an integer such that $- \frac{a_{i_0}}{r_{i_0}} > 
\sum_{i=1}^{n+1} \frac{a_ik_i}{r_i} \ge 0$. 
Then we have 
\[
- \sum_{i=n''+1}^{n+1} \frac{a_i}{r_i} > \sum_{i=1}^{n+1} \frac{a_ik_i}{r_i} 
- \sum_{p=1}^t \frac{a_{i_p}\epsilon_p}{r_{i_p}} \ge 0 
\]
if $\epsilon_p = 0, 1$.
Therefore, $\mathcal{D}$ does not change if we add all the sheaves of the form
$\mathcal{O}_S(\sum_{i=1}^{n+1} \frac{k_i}{r_i}D'_i)$ to the set $\Omega$.

Assume that $t = n-n''-1$. Then $S$ is affine and 
$\mathcal{O}_S(\frac 1{r_{i_0}}D'_{i_0}) \cong \mathcal{O}_S$. 
Therefore, $\mathcal{D}$ does not change if we add the invertible sheaves 
$\mathcal{O}_S(\sum_{i=1}^{n+1} \frac {k_i}{r_i}D'_i)$
for all integers $k_i$ to the set $\Omega$.
Hence any coherent sheaf whose support is contained in $S$ belongs to 
$\mathcal{D}$.

We claim that, for any staratum $S$, 
$\mathcal{D}$ does not change if we add the invertible sheaves 
$\mathcal{O}_S(\sum_{i=1}^{n+1} \frac{k_i}{r_i}D'_i)$
for all integers $k_i$ to the set $\Omega$.
Indeed, we proceed by the descending induction on $t$.
Let $S' = S \cap \frac 1{r_{i_0}}D'_{i_0}$.
Then our claim follows from the following exact sequence
\[
0 \to \mathcal{O}_S(- \frac 1{r_{i_0}}D'_{i_0}) 
\to \mathcal{O}_S \to \mathcal{O}_{S'} \to 0.
\]
Hence any coherent sheaf belongs to $\mathcal{D}$.
\end{proof}

\begin{Lem}\label{span'}
The set of all the invertible sheaves on $\mathcal{X}$
in the range \ref{range2} is a spanning class of 
$D^b(\text{Coh}(\mathcal{X}))$.
\end{Lem}

\begin{proof}
The proof is similar to that of Lemma~\ref{span}.
\end{proof}

\begin{Cor}
Let $X$ and $Y$ be quasi-smooth projective toric varieties, 
let $B$ and $C$ be effective toric $\mathbb{Q}$-divisors
on $X$ and $Y$, respectively, whose coefficients are contained in the set 
$\{1 - 1/r \vert r \in \mathbb{Z}_{>0}\}$, and let $\mathcal{X}$ and 
$\mathcal{Y}$ be smooth Deligne-Mumford stacks attached to the 
pairs $(X, B)$ and $(Y, C)$, respectively. 
Let $f: X -\to Y$ be a toric proper birational map which is log crepant 
in the sense that 
\[
g^*(K_X+B) = h^*(K_Y + C)
\]
for toric proper birational morphisms from a common toric variety 
$g: Z \to X$ and $h: Z \to Y$ such that $f = h \circ g^{-1}$.
Then there is an equivalence of triangulated categories
\[
F: D^b(\text{Coh}(\mathcal{Y})) \to D^b(\text{Coh}(\mathcal{X})).
\]
\end{Cor}

\begin{proof}
By \cite{Matsuki} with some additional argument due to Matsuki, 
$f$ is shown to be decomposed into a sequence of divisorial contractions and
flips which are log crepant.
\end{proof}

\begin{Rem}
(1) We note that the Fourier-Mukai functors 
for both divisorial contractions and
flips are of the same type consisting of pull-backs and push-downs.
Indeed, divisorial contractions and
flips are of the same kind of operations from the 
view point of the Minimal Model Program though they look very different
geometrically.

(2) In the situation of Theorem~\ref{toroidal}~(1), 
suppose that there is a third $\mathbb{Q}$-divisor $D$ on $X$ 
with standard coefficients such that $B \ge C \ge D$.
Let $\mathcal{Z}$ be the stack corresponding to the pair $(X, D)$.
Then there are fully faithful functors 
$F_1: D^b(\text{Coh}(\mathcal{Z})) \to D^b(\text{Coh}(\mathcal{Y}))$, 
$F_2: D^b(\text{Coh}(\mathcal{Y})) \to D^b(\text{Coh}(\mathcal{X}))$ 
and $F_3: D^b(\text{Coh}(\mathcal{Z})) \to D^b(\text{Coh}(\mathcal{X}))$ as 
proved there.
But we have $F_3 \not\cong F_2 \circ F_1$ in general.

(3) Let $X = \mathbb{C}^2$, $\bar B = B_1+B_2$ the union of the 
two coordinate lines, 
$f: Y \to X$ the blowing up at the singular point $B_1 \cap B_2$ of $\bar B$, 
and $\bar C = B'_1+B'_2+C_3$ the union of the strict transforms and the 
exceptional divisor.
Let $\mathcal{X}_n$ and $\mathcal{Y}_n$ be the stacks associated to the 
pairs $(X, \frac 1{2n}(B_1+B_2))$ and $(Y, \frac 1{2n}(B'_1+B'_2) + 
\frac 1n C_3)$, respectively.
Since $f^*(K_X+\frac 1{2n}(B_1+B_2)) = K_Y + \frac 1{2n}(B'_1+B'_2) + 
\frac 1n C_3$, we have equivalences $\Phi_n: D^b(\text{Coh}(\mathcal{X}_n)) 
\to D^b(\text{Coh}(\mathcal{Y}_n))$ for each $n$ by 
Theorem~\ref{toroidal}~(2). 
On the other hand, we have fully faithful functors
$\Psi_{X, nn'}: D^b(\text{Coh}(\mathcal{X}_n)) \to 
D^b(\text{Coh}(\mathcal{X}_{n'}))$ and 
$\Psi_{Y, nn'}: D^b(\text{Coh}(\mathcal{Y}_n)) \to 
D^b(\text{Coh}(\mathcal{Y}_{n'}))$ for $n < n'$ by Theorem~\ref{toroidal}~(1). 
We have to be careful that the following diagram is not commutative
\[
\begin{CD}
D^b(\text{Coh}(\mathcal{X}_n)) @>{\Phi_n}>> D^b(\text{Coh}(\mathcal{Y}_n)) \\
@V{\Psi_{X, nn'}}VV @VV{\Psi_{Y, nn'}}V \\
D^b(\text{Coh}(\mathcal{X}_{n'})) @>{\Phi_{n'}}>> 
D^b(\text{Coh}(\mathcal{Y}_{n'})).
\end{CD}
\]
The reason is that the Serre functors, defined by using 
different invertible sheaves, are not compatible.
\end{Rem}

%%%%%%%%%%%%%%%%%%%%%%%%%%%%%%%%%%%%%%%%%%%%%%%%%%%%%%%%%%%%%%%%%%

We conclude this paper with a remark on the  
non-commutative geometry as in \cite{VDBergh}.
We consider the situation of Theorem~\ref{toroidal}~(3) under the additional 
assumption that $\mu^*(K_X+B)=\nu^*(K_Y+C)$.

Although the set of all invertible sheaves on $\mathcal{Y}$ in the range
(\ref{range}) is infinite, there are only finitely many ismomorphism classes.
Let $P_Y$ be the direct sum of these representatives, and let
$A_Y = \text{Hom}(P_Y, P_Y)$ be the non-commutative ring of endomorphisms.
We denote by $\text{Mod}(A_Y)$ the abelian category of finitely 
generated right $A_Y$-modules.

\begin{Prop}\label{non-comm}
There is an equivalence of triangulated categories:
\[
D^b(\text{Coh}(\mathcal{Y})) \cong D^b(\text{Mod}(A_Y)).
\]
\end{Prop}

\begin{proof}
By the proof of Theorem~\ref{toroidal}~(3),
we have $\text{Hom}^p(P_Y, P_Y) = 0$ for $p > 0$ and the set $\{P_Y\}$ 
is spanning.
We claim that the functors
\[
\begin{split}
&G: D^b(\text{Coh}(\mathcal{Y})) \to D^b(\text{Mod}(A_Y)) \\
&H: D^b(\text{Mod}(A_Y)) \to D^b(\text{Coh}(\mathcal{Y}))
\end{split}
\]
given by $G(a) = R\text{Hom}(P_Y, a)$ and $H(m) = m \otimes_{A_Y}^L P_Y$
are quasi-inverses each other.
Indeed, we have $G(P_Y) \cong A_Y$ and $H(A_Y) \cong P_Y$.
Since $P_Y$ and $A_Y$ respectively span $D^b(\text{Coh}(\mathcal{Y}))$ and 
$D^b(\text{Mod}(A_Y))$, it follows that 
$G$ and $H$ are fully faithful by \cite{Bridgeland}, hence quasi-inverses
each other.
\end{proof}

%%%%%%%%%%%%%%%%%%%%%%%%%%%%%%%%%%%%%%%%%%%%%%%%%%%%%%%%%%%%%%%%%%

Department of Mathematical Sciences, University of Tokyo, 

Komaba, Meguro, Tokyo, 153-8914, Japan 

kawamata@ms.u-tokyo.ac.jp


\begin{thebibliography}{ACFMT}

\bibitem{Bridgeland}
T. Bridgeland.
{\em Equivalences of triangulated categories and Fourier-Mukai transforms}.
math.AG/9809114. 
Bull. London Math. Soc. {\bf 31}(1999), 25--34.

\bibitem{BKR}
T. Bridgeland, A. King and M. Reid.
{\em Mukai implies McKay: the McKay correspondence as an equivalence of 
derived categories}.
math.AG/9908027.
J. Amer. Math. Soc. 14 (2001), 535--554.

\bibitem{Caldararu}
A. Caldararu.
{\em The Mukai pairing, I: The Hochschield structure}.
math.AG/0308079.

\bibitem{Francia}
Y. Kawamata.
{\em Francia's flip and derived categories}.
math.AG/0111041,
in Algebraic Geometry (a volume in Memory of Paolo Francia),
Walter de Gruyter, 2002, 197--215.

\bibitem{DK}
Y. Kawamata.
{\em D-equivalence and K-equivalence}.  
math.AG/0205287, J. Diff. Geom. {\bf 61} (2002), 147--171.

\bibitem{equi}
Y. Kawamata.
{\em Euivalences of derived catgories of sheaves on smooth stacks}.
math.AG/0210439,
to appear in Amer. J. Math.

\bibitem{KMM}
Y. Kawamata, K. Matsuda and K. Matsuki. 
{\em Introduction to the minimal model problem}. 
in Algebraic Geometry Sendai 1985,
Advanced Studies in Pure Math. {\bf 10} (1987), 
Kinokuniya and North-Holland, 283--360. 

\bibitem{Matsuki}
K. Matuski.
{\em Introduction to Mori Program}.
Springer, 2002.

\bibitem{VDBergh}
M. Van den Bergh.
{\em Three-dimensional flops and non-commutative rings}.
math.AG/0207170.

\end{thebibliography}
\end{document}